\documentclass{article}
\usepackage{amssymb}

\makeatletter
\@ifundefined{mathbb}{\let\mathfrak=\frak \let\mathbb=\Bbb}{}
\makeatother

\def\fa#1{{\forall #1}\,}
\def\Sym{\mathop{\rm Sym}\nolimits}
\title{Errata for {\it Automorphisms of first-order structures}}
\author{R.~W.~Kaye and H.~D.~Macpherson (eds),\\ 
Oxford University Press, 1994.}
\date{20th June, 1995}
\begin{document}
\maketitle
\bf

\noindent{Page 78, line 23.} Delete the words `{\rm are called}'.

\noindent{Page 93, line 26.} For `{\rm is brief}' read `{\rm will be brief}'.

\noindent{Page 97, line 4 up.} For $\mathbb N$ read ${\rm Sym}(\mathbb N)$.


\noindent{Page 220.} Replace the first seven lines (from `{\rm The
subgroups $K$ of \ldots}' to `{\rm $\alpha_x = id )\}$.}') with the following
passage.
\medskip

{\rm 
\begin{quote}
The subgroups $K$ of $\prod _{x\in N}\Sym(F)$ which are kernels of finite\\ 
covers of $N$ are exactly the groups of the form 
$$
K^{L}_{H} = \{ \alpha  \in  \prod _{x\in N} L: \fa{x,y}%
  (\alpha _{x}H = \alpha _{y}H) \},
$$
where $L$ is a subgroup of $\Sym(F)$, $H$ is a  normal  subgroup  of $L$\\ 
and $\alpha _{a}$ is the action on $F$ which corresponds to the projection\\ 
of $\alpha \hbox{ on }F_{a} = \{a\}\times F$.
\end{quote}
 
If $K$ is a kernel of some cover then by Ziegler (1992) the groups $L$ and\\
$H$ can be found in the following way: let $a \in  N$, then}
\begin{eqnarray*}
L(K) & = & \{\alpha _{a}: \alpha  \in  K\}, \\
H(K) & = & \{\alpha _{a}: \alpha  \in  K\hbox{\rm\ and }\fa {x \neq a} %
           (\alpha _{x} =id)\}.
\end{eqnarray*}

\noindent{Page 269, line 4 up.} For $frakB$ read $\mathfrak B$.

\noindent{Page 285, line 17.} For $0<p\leqslant q\leqslant k$
read $0\leqslant p < q\leqslant k$.

\noindent{Page 298, line 20.} For {\rm anomolous} read {\rm anomalous}.

\noindent{Page 311, line 23.} For $M \vDash {\rm Th}(\mathbb N)$ read
$M \not\vDash {\rm Th}(\mathbb N)$

\noindent{Page 320, line 4 up.}  For {\it action
of a subgroup of transitive permutation group} read {\it
action of a subgroup of a transitive permutation group}.

\end{document}